\documentclass [12pt]{article}
\newfont{\bbb} {msbm10}
\newcommand{\Bbb}[1]{\mbox{\bbb#1}}
\newcommand{\R}{\Bbb{R}}
\newcommand{\bS}{\Bbb{S}}

\newcommand{\Z}{\Bbb{Z}}

\newcommand{\C}{\Bbb{C}}
\newcommand{\Ha}{\Bbb{H}}

\newcommand{\cT}{{\cal{T}}}
\newcommand{\sm}{\setminus}
\newcommand{\sbs}{\subset}
\newcommand{\ra}{\rightarrow}

\newcommand{\met}{{\cal{MET}}}
\newcommand{\mo}{{\cal{MET}}^{\, sec\, < \, 0}}
\newcommand{\Mo}{{\cal{M}}^{\, sec\, < \, 0}}
\newcommand{\To}{{\cal{T}}^{\, sec\, < \, 0}}

\newcommand{\cD}{{\cal{D}}}

\usepackage[all]{xy}

\begin{document}

\title{  The Moduli Space of Negatively Curved Metrics
of a Hyperbolic Manifold}
\author{F. T. Farrell and P. Ontaneda\thanks{Both authors were
partially supported by NSF grants.}}
\date{}

\maketitle


\noindent {\bf \Large  Section 0. Introduction.}\\

Let $M$ be a smooth closed manifold. 
We will denote the group of all self-diffeomorphisms of $M$, with the smooth topology, by $DIFF(M)$.\\

It is natural to study
the space $\mo (M)$
of all (complete) Riemannian metrics on $M$ with negatively curved sectional
curvatures. It is also natural to study $\Mo (M)$, the {\it moduli space of negatively
curved metrics on $M$}. This space is
defined as the quotient space $\mo (M)/\cD(M)$, where $\cD(M)=\R^+\times DIFF(M)$
acts on $\mo (M)$ by scaling and pushing-forward. We can think of $\Mo(M)$ as
the space  of all {\it negatively curved
geometries} on $M$.  Denote by $\kappa$ the quotient map
$$ \mo (M)\stackrel{\kappa}{\longrightarrow} \Mo (M)$$

In this paper we prove that if the closed hyperbolic manifold $M$ has a ``good enough" closed geodesic $\gamma$ then  
$\pi_k(\Mo(M))$ and $H_k(\Mo(M))$ are non-trivial,
provided one of the following hold for $k$ and $n=dim\, M$:

$${\bf \Big(*\Big)}..........\left\{
\begin{array}{l} {\bf 1.}\,\, \,\, k=0\,\,\, {\mbox{and}} \,\,\, n\geq 10\\
{\bf 2.} \,\,\,\, k=1\,\,\, {\mbox{and}} \,\,\, n\geq 12\\
{\bf 3.} \,\,\,\, k=2p-4,\,\, p>2\,\, {\mbox{prime}},\,\,\, {\mbox{and}} \,\,\, n\geq 3k+8\\
\end{array}\right.$$

In fact we prove more. We show that in these cases the maps 

$$\begin{array}{l}
\pi_k \Big(\mo(M)\Big)\stackrel{\pi_k(\kappa)}{\longrightarrow} \pi_k\Big(\Mo(M)\Big)\\ \\ 
 H_k \Big(\mo(M)\Big)\stackrel{H_k(\kappa)}{\longrightarrow}  H_k\Big(\Mo(M)\Big)
\end{array}$$ 
\vspace{.2in}

\noindent are non-zero.\\

\noindent {\bf Remarks.} 

\noindent {\bf 1.} The concept of $\gamma$ being ``good enough" depends on $k$.

\noindent {\bf 2.} Here $H_k(\, \, )$ denotes reduced homology with integer coefficients.

\noindent {\bf 3.} For a space $X$, by ``$\pi_0(X)$ is non-zero"
we mean $X$ is not path-connected. If $f:X\ra Y$, by ``$\pi_0(f)$ is non-zero" we mean ``$\pi_0(f)$ is non-constant".
Hence ``$\pi_0(X)$ is non-zero", 
``$H_0(X)$ is non-zero" and 
``$X$ is not path-connected" are all equivalent. And
``$\pi_0(f)$ is non-zero",
``$H_0(f)$ is non-zero"
and ``$\pi_0(f)$ is non-constant" are all equivalent.

\noindent {\bf 4.} All homotopy groups are based at (the class of) the given hyperbolic metric.\\

Before we give a detailed statement of this result we need some notation.
Let $M$ be a closed hyperbolic manifold of $dim\,>\, 2$. The (finite) group of isometries of $M$ will be denoted by
$ISO(M)$. If $\gamma$ is an embedded closed geodesic in $M$,
we will denote by $\omega (\gamma)$ the width of its normal geodesic neighborhood and by $o(\gamma)$ the order of
the group of isometries that preserve $\gamma$, that is $$o(\gamma)=\Bigg|\,\Big\{ f\in ISO(M)\, :\, f(\gamma)=\gamma  \Big\}\,  \Bigg|$$

\noindent Also, the length of $\gamma$ will be denoted by $\ell(\gamma)$. Here is the statement of our main result:\\

{\bf Main Theorem.} {\it There is a function $r=r(n,k,\ell)$, defined for $\ell>0$ and $(k, n)$ satisfying (*), such that
the following holds. \\

Let $k$ be a non-negative integer and $M$ a closed hyperbolic manifold, with $(k, dim\, M)$ satisfying (*). Assume $M$ has
has an embedded closed geodesic $\gamma$, with orientable (hence trivial) normal bundle, satisfying the following two conditions:}

\begin{enumerate}
\item[(i)] $\omega(\gamma)\geq r(dim\, M, k, \ell(\gamma)\, )$

\item[(ii)] {\it  For $k=0,1$, $o(\gamma)$ is odd. For $k=2p-4$, $o(\gamma)$ is odd and not divisible by $p$.}
\end{enumerate}

\noindent {\it Then $\pi_k(\kappa)$ and $H_k(\kappa)$ are non-zero.
In particular $\pi_k(\,\Mo(M)\, )$ and $H_k(\,\Mo(M)\, )$ are non-trivial.}\\

We say that an embedded closed geodesic, with orientable normal bundle, is $k$-{\it good} if it satisfies conditions (i), (ii) above. Therefore,
our Main Theorem can be paraphrased in the following way: If a closed hyperbolic manifold has a $k$-good geodesic,
and $(k, dim\, M)$ satisfy (*), then $\pi_k(\kappa)$ and $H_k(\kappa)$ are non-zero.\\

The next result shows that, for non-arithmetic hyperbolic manifolds, the condition ``have a $k$-good geodesic"
is obtained after taking finite covers.\\

\noindent {\bf Theorem 1.} {\it Let $M$ be a closed non-arithmetic hyperbolic manifold and 
$k$ a non-negative integer, with $(k, dim\, M)$ satisfying (*). Then $M$ has a finite sheeted cover 
that has a $k$-good geodesic.}\\

\noindent {\bf Remark.} An important ingredient in the proof of Theorem 1 is Theorem 3.1 in Section 3.
This Theorem states that every closed hyperbolic manifold contains a closed geodesic $\gamma$ with
$o(\gamma)=1$. In fact, the Addendum of Theorem 3.1 states that the number $P^*(t)$ of such geodesics,
of length $\leq t$, grows at least as fast as a function of the form $c\,\frac{e^{ht}}{t}$, for some $h,\,c>0$.\\ 

\noindent {\bf Corollary.} {\it Let $M$ be a closed non-arithmetic hyperbolic manifold and 
$k$ a non-negative integer, with $(k, dim\, M)$ satisfying (*).
Then $M$ has a finite sheeted cover $N$ such that the maps}
$$\begin{array}{l}
\pi_k \Big(\mo(N)\Big)\stackrel{\pi_k(\kappa)}{\longrightarrow} \pi_k\Big(\Mo(N)\Big)\\ \\ 
 H_k \Big(\mo(N)\Big)\stackrel{H_k(\kappa)}{\longrightarrow}  H_k\Big(\Mo(N)\Big)
\end{array}$$ 
\noindent {\it are non-zero. In particular $\pi_k(\, \Mo(N)\, )$ and $H_k(\,\Mo(N)\, )$ are non-trivial.}\\

The statements of the Main Theorem (and its Corollary)  holds also for $\epsilon$-pinched negatively curved metrics:\\

\noindent {\bf Addendum to the Main Theorem.} {\it The statement of the Main Theorem and its Corollary remain true if we:}
\begin{enumerate}
\item[{i.}] {\it replace the 
decoration ``sec $<$ 0" on $\mo(M)$ in the Main Theorem by ``-1-$\epsilon<sec\leq$-1". }

\item[{ii.}] {\it replace the 
decoration ``sec $<$ 0" on both $\mo(M)$ and $\Mo(M)$ in the Corollary by ``-1-$\epsilon<sec\leq$-1".}
\end{enumerate}
\noindent {\it In these cases we have  $r=r(n,k, \ell, \epsilon )$.}\\

We can interpolate between $\mo(M)$ and $\Mo(M)$ the space
$\To(M)$ called the {\it Teichm\"uller space of negatively curved metrics on $M$}.
This space was defined in \cite{FO5} as  the quotient space $\mo (M)/\cD_0(M)$, where $\cD_0M=\R^+\times DIFF_0(M)$
and $DIFF_0(M)$ is the space of all self-diffeomorphisms on $M$ that are homotopic to
the identity $1_M$. The definition of $\To(M)$ given in \cite{FO5}
 was motivated by the classical definition of the Teichm\"uller space of a surface.\\
 
The group $DIFF(M)/DIFF_0(M)$ acts on $\To(M)$, and the quotient space is $\Mo(M)$.
 The group $DIFF(M)/DIFF_0(M)$ can be considered as a subgroup of $Out(\pi_1(M))$. 
If $M$ is hyperbolic we actually have $DIFF(M)/DIFF_0(M)=Out(\pi_1(M))$, which is a finite group
 ( $dim\, M>2$). 
Therefore we have the following diagram of quotient spaces:\\

$$\mo(M)\,\,\, \stackrel{\xi }{\longrightarrow}\,\,\, \To(M)\,\,\, \stackrel{\varsigma }{\longrightarrow}\,\,\, \Mo(M) $$\

\noindent  
where the first arrow is the quotient map induced by the action of $\cD_0(M)$ on 
$\mo(M)$ and he second arrow is the quotient map induced by the action of $DIFF(M)/DIFF_0(M)\sbs Out(\pi_1(M))$ on $\To(M)$.
And we have  $\kappa =\varsigma \xi $.
Here are some comments concerning the topology of $\mo(M)$ and $\To(M)$.

\begin{enumerate}
\item[1.] 
It was proved in \cite{FO6} that $\mo(M)$ is never connected, provided $dim\, M\geq 10$,
unless it is empty. In fact,
it always has infinitely many components. Also the homotopy groups $\pi_k(\, \mo(M)\, )$ are not just non-trivial
but not finitely generated, provided $k$ and $n=dim\, M$ satisfy (*).  But these non-zero elements in $\pi_k(\,\mo(M)\, )$ 
 die in $\To(M)$, that is, they are mapped to zero by $\pi_k(\kappa)$.

\item[2.] It was proved in \cite{FO5} that every closed hyperbolic $n$-manifold has a finite sheeted cover $M$
such that $\pi_k(\, \To(M)\, )$ is non-zero,
provided $n+k\equiv 3$, $\, mod\,\, 4$, and $n$ is sufficiently large, depending on $k$.
Here $M$ depends on $k$ and always $k>0$. 
Interestingly, these non-zero elements are not in the image of
$\pi_k(\xi):\pi_k(\, \mo(M)\, )\ra\pi_k(\, \To(M)\, )$.  We do not know if these non-zero elements
survive in $\Mo(M)$.

\item[3.] In \cite{FO7} a version of our Main Theorem was proved with $\To(M)$ instead of $\Mo(M)$.
That is, for a closed hyperbolic manifold $M$
there are non-zero elements in $\pi_k(\, \mo(M)\, )$ that survive in $\pi_k(\, \To(M)\, )$,
provided  $k$ and $n$ satisfy (*) and $M$ has a closed geodesic with
large enough tubular neighborhood.
Our Main Theorem shows that, assuming $M$ has a $k$-good geodesic, these non-zero elements
can be chosen so that they survive all the way to $\pi_k(\, \Mo(M)\, )$.
\end{enumerate}

The paper is structured as follows. In section 1 we recall and introduce some necessary definitions and give some preliminary
results. In section 2 we prove the Main Theorem. In section 3 we prove Theorem 1. In section 4 we prove Lemma 2.1
which states that the Teichm\"uller space $\cT(M)$ is a separable ANR, hence it has the homotopy type of a 
countable CW-complex. Finally in section 5 we prove Proposition 2.3 which is needed in the proof of the Main Theorem.\\

The authors are grateful to Ross Geoghegan for his useful comments and suggestions.
\vspace{.9in}

\noindent {\bf \Large  Section 1. Preliminaries.}\\

In this Section we review some concepts and establish some notation that will be needed later.\\

\noindent {\bf \large  1.1.  The Wreath product.}\\

Let $m$ be a positive integer and  $G$ a finite group acting on $\underbar{{\it m}}=\{ 1,2,...,m\}$.
Let $X$ be a space and let $G$ act on $X^m=\underbrace{X\times ...\times X}_{m\,\,\, times}$ 
permuting the factors, that is, by
$\phi (x_1, x_2, ..., x_m)=(x_{\phi^{-1}(1)}, x_{\phi^{-1}(2)}, ..., x_{\phi^{-1}(m)})$, $\phi\in G$.
Define the {\it ``wreath'' product}:\,\,  $X\wr G=X^m/G$. \\

\noindent (This definition differs from -but is related to- the definition of {\it wreath product} used throughout the
literature, which is given as a quotient of $G\times X^m$.)\\

We will need the following facts:
\begin{enumerate}
\item[{i.}] Let $f:X\ra Y$. Then the map $f^m=f\times...\times f :X^m\ra Y^m$ descends to a map
$f\wr G:X\wr G\ra Y\wr G$ and the following diagram commutes
$$\begin{array}{ccc}
X^m& \stackrel{f^m}{\longrightarrow}& Y^m\\
\downarrow&&\downarrow\\
X\wr G& \stackrel{f\wr G}{\longrightarrow}& Y\wr G
\end{array}
$$

\item[{ii.}] Let $S_m$ be the symmetric group, i.e. the group of all permutations on $\underbar{{\it m}}$. 
The wreath product $X\wr S_m$ is the $m$-symmetric product $SP^m(X)$. (See Section 1.2.)

\item[{iii.}] If $H$ is a subgroup of $G$ then there is an obvious map $X\wr H\longrightarrow X\wr G$,
where it is understood that the $H$ action on $\underbar{{\it m}}$ is induced from the $G$ action.
Here is a particular case: let $G$ be a finite group and order the elements of $G$: \,$\phi_1,....\phi_m$.
This ordering induces an action of $G$ on $\underbar{{\it m}}$: $\phi(i)=j$ if $\phi\phi_i=\phi_j$.
And $G$ can be considered a subgroup of the symmetric group $S_m$, with compatible actions on
$\underbar{{\it m}}$. Therefore, once an ordering of $G$ is given,
there is a map $X\wr G\ra SP^m(X)$, where $m$ is the order of $G$.

\item[{iv.}] Let $G$ be a finite group of order $m$ and let \,$\phi_1,....\phi_m$ be a given ordering of its 
$m$ elements. Then $G$ induces an action on $\underbar{{\it m}}$  as in item (iii) above.
Suppose now that $G$ has already an action on $X$. Define the map $V:X\longrightarrow X^m$, 
$x\mapsto (\phi_1^{-1} x,...,\phi_m ^{-1}x)$. 
And it can be checked that  $V(\phi x)=\phi \, V(x)$, $x\in X$, $\phi\in G$. Therefore $F$ induces a
map $\nu: X/G\ra X\wr G$. This map is clearly injective and it can be checked that, in fact, it
is an embedding.

\end{enumerate}
\vspace{.4in}

\noindent {\bf \large  1.2.  The Infinite Symmetric product.}\\

We will also need some facts about symmetric products. Recall that for a space $X$ and an integer $k>0$
the $k$th- symmetric product $SP^kX$ of $X$ is defined as $X^k/S_k$, where the symmetric group $S_k$
acts on $X^k$ by permuting the factors. Note that $X=SP^1 (X)$. If we fix a basepoint $x_0\in X$ we get embeddings 
$SP^kX\hookrightarrow SP^{k+1}X$,
$[(x_1,...,x_k)]\mapsto [(x_1,...,x_k, x_0)]$. The infinite symmetric product $SP^\infty X$ is the union of all $SP^kX$.
We will use the following facts:
\begin{enumerate}
\item[{i.}] The Dold-Thom Theorem says that for a connected CW complex $X$ and $x_0\in X$ we have $\pi_k (SP^\infty X,x_0)\cong H_k (X)$ (integer coefficients).

\item[{ii.}] Note that the inclusion $X\hookrightarrow SP^\infty X$ factors through $SP^kX$:
$$X\hookrightarrow  SP^k X \hookrightarrow SP^\infty X$$

\item[{iii.}] For any space $X$, there is a map $SP^\infty\Big( SP^\infty (X)\Big)\ra SP^\infty (X)$,
given by concatenation.
And the composition $SP^\infty (X)\ra SP^\infty\Big( SP^\infty (X)\Big)\ra SP^\infty (X)$ is clearly
the identity. It follows that the map $SP^\infty (X)\ra SP^\infty\Big( SP^\infty (X)\Big)$ is injective at the
homology and homotopy levels.

\item[{iv.}] For a connected CW complex $X$, and $x_0\in X$, the map $h:\pi_k(X,x_0)\ra\pi_k(SP^\infty X,x_0)\cong H_k(X)$, induced by the inclusion $X\hookrightarrow SP^\infty X$
is the Hurewicz map. We will also call this inclusion the {\it Hurewicz map} and denote it also by $h$. Let $o$ be a positive
integer. We define the map `` {\it o times the Hurewicz map}", the map $oh:X\ra SP^\infty (X)$, $x\mapsto [(\underbrace{x,...,x}_{o\,\, times}, *,*,*...)]$.
\end{enumerate}

\noindent {\bf Lemma 1.2.1.} {\it We have $\pi_k(oh)=o\,\pi_k(h)$.}\\

\noindent {\bf Proof.} This follows from the fact that $SP^\infty(X)$ is an H-space with
the operation given by concatenation. In particular $(oh)(x)=o\, h(x)$, where we are using additive notation in
the second term. Since the operation in $\pi_k\Big(SP^\infty (X)\Big)$ coincides with the one given by
(pointwise) concatenation, the result follows. 
\vspace{.4in}

\noindent {\bf 1.3. The Simplicial Quotient.}\\

We will borrow the following notation from \cite{FO7}. For more details see Section 1.3 of \cite{FO7}.
Let $X$ be a space and $S(X)$ be its singular simplicial set.  Write $X^\bullet=|\, S(X)|$
where the bars denote ``geometric realization". There is a canonical map $X^\bullet\ra X$
which is a weak homotopy equivalence. If $f:X\ra Y$ is a map then the simplicial map 
$S(f) :S(X)\ra S(Y)$ defines a map $f^\bullet :X^\bullet\ra Y^\bullet$.\\

Let $G$ be a topological group
acting freely on $X$. Then $S(G)$ is a simplicial group acting simplicially on $S(X)$
and we get a simplicial set $S(X)/S(G)$. We define the {\it simplicial quotient} as $X//G=|S(X)/S(G)|$. 
The map $S(X)\ra S(X)/S(G)$ defines a map $X^\bullet\ra X//G$.

\vspace{.8in}

\noindent {\bf \Large  Section 2. Proof of the Main Theorem.}\\

 Let $M$ be a hyperbolic $n$-manifold and
let $\gamma :\bS^1\ra M$ be an embedded closed geodesic of length $\ell$. 
We will denote the image of $\gamma $ also by $\gamma$.
Assume that the normal bundle of $\gamma$ is orientable, hence trivial. 
Let $\omega(\gamma)>0$ be the width of the normal geodesic tubular
neighborhood of $\gamma$.\\

In \cite{FO7} we constructed, using $\gamma$, a map $\Delta_\gamma:\To(M)^\bullet\ra  TOP(\bS^1\times\bS^{n-2}) \,//L$
such that the composition (see Section 1.3 for notation)
$$\mo(M)^\bullet\stackrel{\xi^\bullet}{\longrightarrow}\To(M)^\bullet\stackrel{\Delta_\gamma}{\longrightarrow} TOP(\bS^1\times\bS^{n-2}) \,//L
$$

\noindent is non-zero at the the $\pi_k$ and $H_k$ level, provided $\omega(\gamma)$ is sufficiently large (how large depending on 
$n$, $k$ and $\ell$), where $n$ and $k$ satisfy condition (*) in the introduction.\\

\noindent {\bf Remark.} Here $TOP(\bS^1\times\bS^{n-2})$ is the group of all self-homeomorphisms of $\bS^1\times\bS^{n-2}$ and,
as defined in the Introduction of \cite{FO7}, $L$ is the subgroup of $TOP(\bS^1\times\bS^{n-2})$
of all ``orthogonal" self-homeomorphisms of 
$\bS^1\times\bS^{n-2}$, that is $f:\bS^1\times\bS^{n-2}\ra\bS^1\times\bS^{n-2}$
belongs to $L$
if $f(z,u)=(e^{i\theta}z, A(z)u)$, for some $e^{i\theta}\in \bS^1$, and $A:\bS^1\ra SO(n-1)$.\\

In fact, it can readily be checked from Section 2.3 of \cite{FO7} (see also Theorem C of \cite{FO7}) that we have
more: the composition of maps going from the upper left group to the lower right group in the following
commutative diagram of groups 

$${\small \begin{array}{ccccc}
\pi_k\Big( \mo(M)^\bullet \Big) & \longrightarrow & \pi_k\Big(\To(M)^\bullet\Big) &\longrightarrow & \pi_k\Big( TOP(\bS^1\times\bS^{n-2}) \,//L\Big) \\
\downarrow&& \downarrow&&\downarrow\\
H_k\Big( \mo(M)^\bullet \Big) & \longrightarrow & H_k\Big(\To(M)^\bullet\Big) &\longrightarrow & H_k\Big( TOP(\bS^1\times\bS^{n-2}) \,//L\Big) 
\end{array} }
$$\\

\noindent is non-zero  in the cases mentioned above, that is, provided $\omega(\gamma)$ is sufficiently large, and $n$ and $k$ satisfy condition (*).
Here the vertical arrows are all Hurewicz maps.\\

\noindent {\bf Remark.} In fact, we can choose elements in $\pi_k\Big(\mo(M)^\bullet\Big)$ such that the images in 
$H_k\Big( TOP(\bS^1\times\bS^{n-2}) \,//L\Big)$ have: order 2, for $k=0,\,1$, and
order $p$, for $k=2p-4$. This will be used in Lemma 2.2 and later.\\

Now, we want a map defined on $\To(M)$ and not just on $\To(M)^\bullet$. For this we need the following Lemma.\\

\noindent {\bf Lemma 2.1.} {\it The Teichm\"uller space $\To(M)$ has the homotopy type of a countable CW-complex.}\\

The proof of the Lemma is given in Section 4.\\

It follows from the Lemma above that the natural map $\To(M)^\bullet\ra\To(M)$
(always a weak homotopy equivalence) has a homotopy inverse $\eta:\To(M)\ra\To(M)^\bullet$.
Define $\delta=\delta_\gamma=\Delta_\gamma \,\eta$. Then the composition
$\mo(M)\stackrel{\xi}{\longrightarrow}\To(M)\stackrel{\delta}{\longrightarrow} TOP(\bS^1\times\bS^{n-2}) \,//L$
\,\, is non-zero at the the $\pi_k$ and $H_k$ level,
provided $\omega(\gamma)$ is sufficiently large and  $n$ and $k$ satisfy condition (*).
And, as mentioned before, we have more. We can drop the upper dots in the diagram of groups above:
the composition of maps going from the upper left group to the lower right group in the following
commutative diagram of groups 

$${\small \begin{array}{ccccc}
\pi_k\Big( \mo(M) \Big) & \longrightarrow & \pi_k\Big(\To(M)\Big) &\longrightarrow & \pi_k\Big( TOP(\bS^1\times\bS^{n-2}) \,//L\Big) \\
\downarrow&& \downarrow&&\downarrow\\
H_k\Big( \mo(M) \Big) & \longrightarrow & H_k\Big(\To(M)\Big) &\longrightarrow & H_k\Big( TOP(\bS^1\times\bS^{n-2}) \,//L\Big) 
\end{array} }
$$\\

\noindent is non-zero, provided $\omega(\gamma)$ is sufficiently large, 
and $n$ and $k$ satisfy condition (*). \\

And the claim made in the remark before Lemma 2.1 still holds.
Consider now the Hurewicz maps $h$ and $h'$ (see Section 1.2):
$$ {\scriptsize \begin{array}{ccc} 
TOP(\bS^1\times\bS^{n-2}) \,//L&\stackrel{h}{\longrightarrow}& SP^\infty\Big( \, TOP(\bS^1\times\bS^{n-2}) \,//L\, \Big)\\ \\
SP^\infty\Big(TOP(\bS^1\times\bS^{n-2}) \,//L\Big)&\stackrel{h'}{\longrightarrow}& 
SP^\infty\Bigg(SP^\infty\Big( \, TOP(\bS^1\times\bS^{n-2}) \,//L\, \Big)\,\Bigg)\end{array}}
$$

\noindent {\bf Lemma 2.2.} {\it The composition map}
$$ {\scriptsize \begin{array}{ccccccc} 
\mo(M)&\stackrel{\xi}{\longrightarrow}&\To(M)&\stackrel{\delta}{\longrightarrow}& 
 TOP(\bS^1\times\bS^{n-2}) \,//L&
\stackrel{h}{\longrightarrow} &SP^\infty\Big( \, TOP(\bS^1\times\bS^{n-2}) \,//L\, \Big)\\ \\
&&&&&&\downarrow h'\\ \\
&&&&&& SP^\infty\Bigg(SP^\infty\Big( \, TOP(\bS^1\times\bS^{n-2}) \,//L\, \Big)\,\Bigg)
\end{array} }
$$
\noindent {\it is non-zero at the the $\pi_k$ level,
provided $\omega(\gamma)$ is sufficiently large and  $n$ and $k$ satisfy condition (*).
In fact, we can choose elements in $\pi_k\Big(\mo(M)\Big)$ such that the images in 
$\pi_k\Bigg(SP^\infty\Big( \, TOP(\bS^1\times\bS^{n-2}) \,//L\,\Big) \Bigg)$ have:
order 2 for $k=0,\,1,$ and order $p$, for $k=2p-4$.}\\

\noindent {\bf Proof.}  The map (see Section 1.2)
$${\scriptsize \begin{array}{c} \pi_k\Big( TOP(\bS^1\times\bS^{n-2}) \,//L\Big)
\longrightarrow \pi_k\Bigg( SP^\infty\Big( \, TOP(\bS^1\times\bS^{n-2}) \,//L\,\Big) \Bigg)=
H_k\Big( TOP(\bS^1\times\bS^{n-2}) \,//L\Big)\end{array} }$$
\noindent is just the Hurewicz map. Hence Lemma 2.2  follows from
the assertion made about the second diagram of groups above, and the 
fact that $H_k(X)\ra H_k\Big( SP^\infty (X)\Big)$ is always injective (see item (iii) in Section 1.2). This proves Lemma 2.2.\\

Let $G$ be the isometry group of $M$. The action of $G$ on $\mo(M)$ descends to an action on $\To(M)$, and $\To(M)/G=\Mo(M)$.
Also, since every isometry preserves the given hyperbolic metric (by definition), the action of $G$ preserves the (chosen)
base-points of $\mo(M)$ and $\To(M)$, inducing actions on the corresponding homotopy groups.\\

Let $G_\gamma$ be the subgroup of isometries preserving $\gamma$, and
$G_\gamma^+$ the subgroup of $G_\gamma$ of all isometries that in addition preserve an orientation of $\gamma$,
and preserve an orientation of the normal bundle of $\gamma$.
To simplify our notation, denote by $J$ the composition 

$${\small \begin{array}{c}
\mo(M) \, \stackrel{\xi}{\longrightarrow} \, \To(M) \,\stackrel{\delta}{\longrightarrow}  \, TOP(\bS^1\times\bS^{n-2}) \,//L
\end{array} }
$$

In what follows $\alpha:\bS^k\ra\mo(M)$ will represent an element $[\alpha]\in\pi_k\Big(\mo(M)\Big)$ such that $\pi_k(h'hJ)([\alpha])\neq 0$,
that is, $h'hJ\alpha$ is not null-homotopic. Moreover we will assume that $\pi_k(hJ)([\alpha])$ has order 2, for $k=0,1\,$ and
order $p$, for $k=2p-4$.
The existence of such an $\alpha$ is granted by Lemma 2.2.\\

\noindent {\bf Proposition 2.3.} {\it We can chose $\alpha$ so that it satisfies the following
two conditions:}
\begin{enumerate}
\item[{\bf (a)}] {\it If $f\notin G_\gamma$ then  $Jf\alpha$ is null-homotopic.}

\item[{\bf (b)}] {\it If $f\in G_\gamma^+$ then  $Jf\alpha\simeq J\alpha$, provided $o(\gamma)$ is odd.}
\end{enumerate}

\noindent {\bf Remark.} Note that $f\notin G_\gamma$ if and only if $f^{-1}\notin G_\gamma$. Similarly for
$G_\gamma^+$.\\

The proof of this Proposition requires a detailed use of the techniques and constructions given in \cite{FO7}.
The proof is given in Section 5.\\

Recall that in the introduction we denoted the order of $G_\gamma$ by $o=o(\gamma)$. In our Main Theorem
we demand $o$ to be odd, so, in what remains of this Section, we assume this. This implies that
there are no isometries preserving $\gamma$ that  reverse the orientation of $\gamma$, 
or that reverse the orientation ``transverse" to $\gamma$ (otherwise the order would be even). 
It follows that $G_\gamma=G_\gamma^+$.\\

Let $m$ be the order of $G$. Order the elements of $G$: $f_1=1_M, f_2,...,f_o,...,f_m$
such that the first $o$ elements correspond to the elements of $G_\gamma$.
Let $\mu:\To(M)\ra \Big(\To(M)\Big)^m$ be defined by $[g]\mapsto \Big(f_1^{-1}[g], f_2^{-1}[g],...,f_m^{-1}[g]\Big)$,
for $g\in\mo(M)$.
Also, let $\zeta: TOP(\bS^1\times\bS^{n-2}) \,//L\ra \Big( TOP(\bS^1\times\bS^{n-2}) \,//L\Big)^m$
be defined by $x\mapsto (x,x,...,x, *,...,*)$, for $x\in TOP(\bS^1\times\bS^{n-2}) \,//L$. That is,
the first $o$ components of $\zeta(x)$ are equal to $x$, and the remaining ones are equal to 
the base-point $*$ of $ TOP(\bS^1\times\bS^{n-2}) \,//L$. We have the following diagram

$${\small \begin{array}{ccccccc}
\bS^k& \stackrel{\alpha}{\longrightarrow}&\mo(M)&\stackrel{\xi}{\longrightarrow}&\To(M)
&\stackrel{\delta}{\longrightarrow}& TOP(\bS^1\times\bS^{n-2}) \,//L\\ \\
&&&&\mu\downarrow&&\downarrow\zeta\\ \\
&&&&\Big(\To(M)\Big)^m
&\stackrel{\delta^m}{\longrightarrow}& \Big(TOP(\bS^1\times\bS^{n-2}) \,//L\Big)^m
\end{array}}
$$\\

\noindent{\bf Lemma 2.4.} {\it We can choose $\alpha$  
such that $\delta^m\mu\,\xi\,\alpha\,\simeq\,\zeta\,\delta\,\xi\,\alpha$.}\\

\noindent {\bf Proof.} It follows directly from Proposition 2.3 and the fact that we are
assuming $G_\gamma=G_\gamma^+$ (see also the Remark after the statement of Proposition 2.3). This proves the Lemma.\\

This last diagram embeds in the following bigger one:

$${\small \begin{array}{ccccccc}
\bS^k& \stackrel{\alpha}{\longrightarrow}&\mo(M)&\stackrel{\xi}{\longrightarrow}&\To(M)
&\stackrel{\delta}{\longrightarrow}& TOP(\bS^1\times\bS^{n-2}) \,//L\\ \\
&&&^\varsigma\swarrow&\mu\downarrow&&\downarrow\zeta\\ \\
&&\Mo(M)&&\Big(\To(M)\Big)^m
&\stackrel{\delta^m}{\longrightarrow}& \Big(TOP(\bS^1\times\bS^{n-2}) \,//L\Big)^m\\ \\
&&&\nu\searrow&q\downarrow&&\downarrow q'\\ \\
&&&&\To(M) \wr G
&\stackrel{\delta\wr G}{\longrightarrow}& TOP(\bS^1\times\bS^{n-2}) \,//L  \wr G\\ \\
&&&&&&\downarrow\sigma\\ \\
&&&&&&SP^\infty\Big(  TOP(\bS^1\times\bS^{n-2}) \,//L \Big)
\end{array}}
$$\\

\noindent where $q$ and $q'$ are quotient maps (see Section 1.1), and $\sigma$ is the composition
(see Section 1.1 and 1.2)

$${\scriptsize \begin{array}{c}
TOP(\bS^1\times\bS^{n-2}) \,//L  \wr G\longrightarrow 
 SP^m\Big( TOP(\bS^1\times\bS^{n-2}) \,//L \Big)\longrightarrow SP^\infty\Big(  TOP(\bS^1\times\bS^{n-2}) \,//L \Big)
\end{array}}
$$\\

\noindent Also, $\nu$ is the natural embedding and the commutative left triangle 
(which can also be viewed as a square) follows from item (iv) in Section 1.1. Note that the composition
$\sigma\,q'\,\zeta$ of all right-hand vertical arrows is nothing but $oh$, i.e, $o$ times the Hurewicz map $h$.
Write $\theta=\sigma\,(\delta\wr G)$.
Since the lower square commutes (see Section 1.2) Lemma 2.4 implies the following Corollary\\

\noindent{\bf Corollary 2.5.} {\it We can choose $\alpha$ 
such that $(oh)\delta\,\xi\,\alpha\,\simeq\,\theta\,q\,\mu\,\xi\,\alpha$.}\\

Write $\varrho=\sigma\,(\delta\wr G)\,\nu$. Then, from the big diagram above we can extract the following
sub-diagram, which is our main diagram:

$${\small \begin{array}{ccccc}
\mo(M)&\stackrel{\xi}{\longrightarrow}& \To(M)&\stackrel{\varsigma}{\longrightarrow}&\Mo(M)\\ \\
&&\delta\downarrow&&\downarrow\varrho\\ \\ 
&& TOP(\bS^1\times\bS^{n-2}) \,//L &\stackrel{oh}{\longrightarrow}&SP^\infty\Big(  TOP(\bS^1\times\bS^{n-2}) \,//L \Big)
\end{array} }$$\\


Now, this diagram is not commutative, but Corollary 2.5 says that we can find $\alpha:\bS^k\ra\mo(M)$ such
that $\pi_k\Big( (oh)\delta\,\xi\Big)([\alpha])\,=\,\pi_k\Big(\varrho\,\varsigma\,\xi\Big)([\alpha])$.\\

But $\pi_k\Big( (oh)\delta\,\xi\Big)([\alpha])\,=\,o\,\pi_k\Big( h\,\delta\,\xi\Big)([\alpha])\,=\,
o\pi_k\Big( hJ\Big)([\alpha])$, and recall that the order of $\pi_k\Big( hJ\Big)([\alpha])$ is
2, for $k=0,\, 1$ and $p$, for $k=2p-4$. Therefore, if $o$ is odd (and not divisible by $p$, for $k=2p-4$)
it follows that $\pi_k\Big( (oh)\delta\,\xi\Big)([\alpha])$ is not zero. Consequently $\pi_k(\kappa)([\alpha])=
\pi_k(\varsigma\,\xi)([\alpha])$ is not zero. This proves the Main Theorem for the homotopy group case.\\

To prove the homology case, apply the $SP^\infty (\,.\, )$ functor to the main diagram above.
And there are corresponding maps from the main diagram to the $SP^\infty$ version, and everything commutes. But we have:

\begin{enumerate}
\item[{\bf 1.}] $\pi_k\Big( h'hJ\Big)([\alpha])$ is not zero (see Lemma 2.2).

\item[{\bf 2.}] the order of $\pi_k\Big( hJ\Big)([\alpha])$ is
2, for $k=0,\, 1$ and $p$, for $k=2p-4$.
\end{enumerate}

Therefore the order of $\pi_k\Big( h'hJ\Big)([\alpha])$ is
2, for $k=0,\, 1$ and $p$, for $k=2p-4$. Hence we can apply the same argument to the $SP^\infty$ version of 
the main diagram

$${\scriptsize \begin{array}{ccccc}
SP^\infty\Big(\mo(M)\Big)&\longrightarrow& SP^\infty\Big(\To(M)\Big)
&\longrightarrow&
SP^\infty\Big(\Mo(M)\Big)\\ \\
&&\downarrow&&\downarrow \\ \\ 
&&SP^\infty\Big( TOP(\bS^1\times\bS^{n-2}) \,//L\Big) &\longrightarrow&
SP^\infty\Bigg(SP^\infty\Big(  TOP(\bS^1\times\bS^{n-2}) \,//L \Big)\Bigg)
\end{array} }$$

\noindent and show that the map 

$${\scriptsize \begin{array}{c}\pi_k\Bigg( SP^\infty \Big(\mo(M)\Big)\Bigg)\ra
\pi_k\Bigg( SP^\infty \Big(\To(M)\Big)\Bigg) \ra \pi_k\Bigg( SP^\infty \Big(\Mo(M)\Big)\Bigg)
\end{array} }$$\\

\noindent is not zero. Hence the homology case follows. This proves our Main Theorem.
\vspace{.8in}

\noindent {\bf \Large  Section 3. Proof of Theorem 1.}\\

In this section $M$ will denote a closed connected (real) hyperbolic $n$-manifold, $n>2$,
and $G=ISO(M)$. Also in this section, by a {\it geometric closed geodesic} of $M$  we mean the image $\gamma$ of a totally-geodesic immersion $\gamma:\bS^1\ra M$ (we use the same symbol for the map and its image).
And a by a {\it simple closed geodesic} we mean a totally-geodesic connected dimension-one submanifold of $M$.
Note that geometric (and simple) closed geodesics have a well-defined length.
Here is the main ingredient of the proof:\\

\noindent {\bf Theorem  3.1.} {\it $M$ contains a
geometric closed geodesic $\gamma$ such that $o(\gamma)=1$, i.e. $G_\gamma$ is trivial.}\\

The method used to prove Theorem 3.1 yields in fact a stronger result about the growth rate of the
number of geometric closed geodesics $\gamma$ in $M$, of length $\leq t$,
with $o(\gamma)=1$.  Let's denote this number by $P^*(t)$.\\

\noindent {\bf Addendum to Theorem  3.1.} {\it There is a constant $c>0$ such that, for $t>1$,}
$$ 
c\,\, \frac{e^{\, ht}}{t}\,\, \, \leq \,\,\, P^*(t)
$$
\\

\noindent Here the positive constant $h=h_n$ depends only on $n=dim\, \Ha^n$ and is the asymptotic growth rate of the volume of balls of radius $t$ in $\Ha^n$, that is:
$$
h\,\,=\,\, h_n\,\, =\,\, lim_{t\ra\infty}\,\,\frac{ln\Big({\mbox{{\bf Vol}}}_{\,\Ha^n}( {\mbox{ball of radius }}\, t \,)\Big)}{t}
$$

\noindent {\bf Remark.} The number of geometric closed geodesics in $M$ of length $\leq t$
is usually denoted by $P_{reg}(t)$, or simply by $P(t)$.
In section 5.3 of \cite{K} it is shown that there are positive constants $a$, $b$ (depending on $M$) 
such that, for $t>1$,
$$
a\,\, \frac{e^{\, ht}}{t}\,\, \, \leq \,\,\, P(t)\,\,\, \leq\,\,\, b\,\,\frac{e^{\, ht}}{t}
$$
\\

Before we prove Theorem 3.1 and its Addendum, we show how it implies Theorem 1. In fact we show how Theorem 3.1
implies the following result which is a bit more general than Theorem 1.\\

\noindent {\bf Theorem 3.2.} {\it Let $M$ be non-arithmetic. Then there is an orientable finite sheeted cover $N$
of $M$ together with a simple closed geodesic $\gamma$ in $N$ such that the following holds.
Given $r>0$, there are a finite sheeted cover $P\ra N$ and a lifting  $\gamma'$ of $\gamma$ to $P$
such that $\omega (\gamma')>r$ and $o(\gamma')=1$.}\\

\noindent {\bf Proof of Theorem 3.2.}
We have that $M=\Ha^n/\Pi$, for some discrete subgroup $\Pi$ of $ISO(\Ha^n)$.  And since
$M$ is non-arithmetic, there exists a unique maximal lattice $\Gamma$ in $ISO(\Ha^n)$
containing $\Pi$, namely $\Gamma$ is the commensurator of $\Pi$, and $\Pi$ has finite index in $\Gamma$.
Let $\Pi_0$ be a finite index
normal subgroup of $\Gamma$ contained in $\Pi$ and set $N=\Ha^n/\Pi_0$. 
We can arrange that $N$ is orientable by first replacing $M$ by its oriented cover, if necessary.
Apply Theorem 3.1 to
$N$ to obtain a geometric closed geodesic $\gamma$ in $N$ such that $ISO(N)_\gamma$ is trivial,
that is, no isometry of $N$ preserves $\gamma$.\\

Let $r>0$. We can apply Corollary 3.3 of \cite{FJ6} to obtain a finite sheeted cover $p:P\ra N$
and a lifting $\gamma'$ of $\gamma$ to $P$ which is a simple closed geodesic and such that $\omega (\gamma')>r$. To finish
the proof of the Theorem we prove that $o(\gamma')=1$.\\

Let $f\in ISO(P)_{\gamma'}$.
Let $L$ be a line in $\Ha^n$ covering $\gamma'$ under the covering projection
$\Ha^n\ra\Ha^n/\pi_1(P)=P$ and let $\hat{f}\in ISO(\Ha^n)$ be a lift of $f$ which preserves $L$.
Since $\Gamma$ is the commensurator of $\Pi$ and $\pi_1(P)$ has finite index in $\Pi$,
we have that $\hat{f}\in \Gamma$. But $\Pi_0$ is normal in $\Gamma$, hence
$\hat{f}$ induces an isometry $\bar{f}:N\ra N$ such that the diagram below commutes:
$$\begin{array}{ccc}
P&\stackrel{f}{\longrightarrow}& P\\
p\downarrow\,\,\, &&\,\,\,\downarrow p\\
N& \stackrel{\bar{f}}{\longrightarrow}& N
\end{array}$$

\noindent But $\bar{f}$ leaves $\gamma$ invariant, hence it is the identity $1_N$.
Therefore $f$ is a covering transformation that leaves $\gamma'$ invariant. 
Since a covering transformation (other than the identity) has no fixed points
and for $q=p|_{\gamma'}:\gamma'\ra\gamma$ there is a point in $\gamma$ whose pre-image under $q$ has
cardinality one, we conclude that $f=1_P$.
This proves Theorem 3.2.\\

To finish the proof of Theorem 1 (and Theorem 3.2) we have to prove Theorem 3.1. For this
we shall use the following definitions. \\

\begin{enumerate}
\item[{\bf 1.}] As mentioned above, the number of geometric closed geodesics in $M$, of length $\leq t$,
is denoted by $P(t)$. And for a non-empty open subset $U$ of $M$,
$P_U(t)$ will denote the number of geometric closed geodesics in $M$ that meet $U$.
\end{enumerate}

In section 5.3 of \cite{K} it is also shown that there as a positive constant  $a_U$ (depending on $U$) such that, for $t>1$,\\

$
(*)\hspace{1.5in}
a_U\,\, \frac{e^{\, ht}}{t}\,\, \, \leq \,\,\, P_U(t)\,\,\, \leq\,\,\, b\,\,\frac{e^{\, ht}}{t}
$\\

\begin{enumerate}
\item[{\bf 2.}] We will denote by $P^T(t)$ the number of geometric closed geodesics $\gamma$ in $M$, of length $\leq t$,
that {\it have a translation}, that is, there is $f\in G_{\gamma}$ such that $f|_{\gamma}:\gamma\ra\gamma$ is a non-trivial
translation, i.e. there is $z_0\in\bS^1\sm\{ 1\}$ such that $f|_\gamma(\gamma(z))=\gamma(z_0\,z)$, for all $z\in\bS^1\subset\C$.
The number $P^T_U(t)$ is the number of such geodesics that meet $U$.

\item[{\bf 3.}] Define $G_\gamma'=\{ f|_{\gamma}\,\,  ;\,\, f\in G_\gamma\}$. We say that $\gamma$ is a {\it fixed}
geodesic if $G_\gamma\ra G_\gamma'$ is not an isomorphism. Equivalently, $\gamma$ is contained in the fixed
point set of some non-trivial isometry $f\in G$. If this is not the case we say that $\gamma$ is a {\it non-fixed} geodesic.
Let's denote by $P^{fixed}(t)$
the number of fixed geometric closed geodesics in $M$, of length $\leq t$.

\item[{\bf 4.}] And we say that a non-fixed geodesic $\gamma$ {\it has no translations}  if
$f\in G_{\gamma}\sm\{ 1_M\}$ implies that $f|_{\gamma}:\gamma\ra\gamma$
reverses orientation, that is 
there is $z_0\in\bS^1\sm\{ 1\}$ such that $f|_\gamma(\gamma(z))=\gamma(z_0\,{\bar{z}})$, for all $z\in\bS^1\subset\C$.
Hence, if $\gamma$ has no translations then either $o(\gamma)=1$ or 
$G_\gamma'\cong G_\gamma$ is isomorphic to $\Z_2$.

\item[{\bf 5.}] Given $f\in G\sm\{ 1_M\}$ 
we define $P^{f}(t)$
as being the number of non-fixed geometric closed geodesics $\gamma$ in $M$, of length $\leq t$,
with no translations, for which $f|_{\gamma}$ is not trivial.
Note that for such geodesics $G_\gamma=\{ 1_M, f\}$ and $f|_\gamma$ reverses orientation.

\item[{\bf 6.}] And recall that we are denoting by $P^*(t)$ the number of geometric closed geodesics 
$\gamma$ in $M$, of length $\leq t$, with $o(\gamma)=1$.
\end{enumerate}

\noindent {\bf Lemma 3.3.} {\it For any non-empty open subset $U$ of $M$ we have, for $t>1$,}
$$
\Bigg( a_U\,\, \frac{e^{\, ht}}{t}\,\, -\,\, P^{fixed}(t)\,\, -\,\, P_U^T(t)\,\, -\,\, \sum_{f\in G\sm\{ 1_M\}}P^f(t)\Bigg)
\,\,\,\leq\,\,\, P^*(t)
$$

\noindent {\bf Proof.} It follows from the definitions above and inequality (*).\\

\noindent {\bf Lemma 3.4.} 
{\it Write $h'=h_{(n-1)}$. There is a constant $c>0$ such that, for $t>1$,}
$$
 P^{fixed}(t)\,\,\, \leq\,\,\, c\,\,\frac{e^{\, h't}}{t}
$$

\noindent {\bf Proof.} The Lemma follows from the formula given in the Remark above, together with the
fact that the fixed-point set of an isometry is a finite disjoint union of closed submanifolds of $M$,
all of dimensions at most $(n-1)$.\\

The following two propositions prove that
$P_U^T(t)$ and $P^{f}(t)$ grow at most as fast as $c\, e^{(2/3)ht}$ does.\\

\noindent {\bf Proposition 3.5. } {\it There is a constant $c>0$ and a non-empty open subset $U$ of $M$, such that, for $t>1$,}
$$
 P^T_U(t)\,\,\, \leq\,\,\, c\,\,\frac{e^{\, (2/3)ht}}{t}
$$
\\

\noindent {\bf Proposition 3.6. } {\it Let $f\in G\sm\{ 1_M\}$. 
Then there is a constant $c>0$ such that, for $t>1$,}
$$
 P^{f}(t)\,\,\, \leq\,\,\, c\,\,\frac{e^{\, (2/3)ht}}{t}
$$
\\

\noindent Because of Lemmas 3.3 and 3.4, and the fact that $0<h'<h$,
these last two Propositions imply Theorem 3.1 and its Addendum.
\\

\noindent {\bf Remark.} The term $(2/3)$ in both formulas in the last two Propositions above
can be improved to $(1/2)+\epsilon$, for any $\epsilon >0$. But then, the constant $c$ may change.\\

In what follows 
we shall denote by {\bf B}$_{\Ha^n}(x,a)$ the closed ball in $\Ha^n$ of radius $a$, centered at $x\in\Ha^n$.\\

\noindent {\bf Proof of Proposition 3.5.} We denote by $p:\Ha^n\ra M=\Ha^n/\Pi$ the covering projection.
Let $U$ be a small  open ball of radius $s>0$ in $M$ consisting of
regular points of the $G$-action on $M$, i.e. if $f\in G\sm\{ 1_M\}$, then $U\cap fU=\emptyset$.
Fix a lift of $U$ to $\Ha^n$, and denote it also by $U$. \\

Define $R_U^T(t)$ as the set of all
oriented simple closed geodesics $\gamma$ in $M$, of length $\leq t$, that have a translation and meet $U$.
By definition we have that $P_U^T(t)$ is the cardinality of $R_U^T(t)$.
As mentioned in item 2 above, for each $\gamma\in R_U^T(t)$ there is an isometry $f\in G$
such that $f|_\gamma$ is a non-trivial translation. For each such $\gamma$ choose one of these
isometries and call it $f_\gamma$. Also, for each $\gamma\in R_U^T(t)$ choose a component $L_\gamma$ of
$p^{-1}(\gamma)$ (i.e. a line) that meets $U\sbs\Ha^n$, and choose the lifting 
$\hat{f}_\gamma$ of $f_\gamma$ to $\Ha^n$ that preserves $L_\gamma$. Then $\hat{f}_\gamma|_{L_\gamma}$
is a translation, and can be chosen to translate by a distance $\leq t/2$.
Let $C\sbs\Ha^n$ compact such that $p(C)=M$ and $U\sbs C$.  Write $U_\gamma=\hat{f}_\gamma (U)$.\\

\noindent {\bf Claim.} {\it If $\gamma\neq \gamma'$ then $U_{\gamma}$ and $U_{\gamma'}$ are disjoint.}\\

\noindent {\bf Proof of the Claim.} If $f_\gamma\neq f_{\gamma'}$ we are done because $U_{\gamma}$ and $U_{\gamma'}$
map, via $p$, to $f_{\gamma}(U)$ and $f_{\gamma'}(U)$, which are disjoint. If $f_\gamma = f_{\gamma'}$
then both $\hat{f}_\gamma$ and  $\hat{f}_{\gamma'}$ are liftings of the same isometry. Hence
$\hat{f}_{\gamma'}=\pi\hat{f}_\gamma$, for some $\pi\in\Pi$. 
But $\gamma\neq\gamma'$ implies $L_\gamma\neq L_{\gamma'}$, hence $\hat{f}_\gamma\neq\hat{f}_{\gamma'}$. Therefore $\pi\neq 1_M$.
Since $\Pi$ acts freely, the Claim follows.\\

Choose a point $x_0\in C$ and let $D$ be the diameter of $C$.
Since $\hat{f}_\gamma|_{L_\gamma}$
is a translation by a distance $\leq t/2$, for each $\gamma\in R_U^T(t)$, we have that
$$
U_{\gamma}\sbs {\mbox{\bf B}}_{\Ha^n}\Big( x_0, (t/2)+s+2D\Big)
$$

\noindent Therefore 

$$
\bigcup_{\gamma\in R_U^T(t)} U_{\gamma}\sbs {\mbox{\bf B}}_{\Ha^n}\Big( x_0, (t/2)+s+2D\Big)
$$

\noindent This, together with the Claim imply

$$
\nu\, P^T_U(t)=\sum_{\gamma\in R_U^T(t)} {\mbox{{\bf Vol}}}\Big(U_{\gamma}\Big)\,\,\, \leq\,\,\, 
{\mbox{{\bf Vol}}}\Bigg({\mbox{\bf B}}_{\Ha^n}\Big( x_0, (t/2)+s+2D\Big)\Bigg)
$$

\noindent where $\nu>0$ is the volume of a ball of radius $s$ in $\Ha^n$. Hence, for $t$ large we have

$$
 P_U^T(t)\,\,\, \leq\,\,\,\frac{1}{\nu}\,\,
{\mbox{{\bf Vol}}}\Bigg({\mbox{\bf B}}_{\Ha^n}\Big( x_0, (8/15)t\Big)\Bigg)
$$

\noindent Taking $ln$ and dividing by $t$ we get

$$
\frac{ln \Big(P_U^T(t)\Big)}{t}\,\,\,\leq\,\,\, -\frac{ln(\nu)}{t}\,\,+\,\, 
(8/15)\,\,\Bigg[\frac{ln \Bigg({\mbox{{\bf Vol}}}\Big({\mbox{\bf B}}_{\Ha^n}\Big( x_0, (8/15)t\Big)\Big)\Bigg)}{(8/15)t}\Bigg]
$$

\noindent For $t$ large we have that the term in the big brackets is less than $(5/4)h$, and we get

$$
\frac{ln \Big(P_U^T(t)\Big)}{t}\,\,\,\leq\,\,\, -\frac{ln(\nu)}{t}\,\,+\,\, (2/3)h
$$

\noindent Multiplying by $t$ and taking the exponential we get

$$
 P_U^T(t)\,\,\, \leq\,\,\, c\,\,\frac{e^{\, (2/3)ht}}{t}
$$

\noindent for $t$ large enough
and $c=\frac{1}{\nu}$. By replacing the constant $c>0$ by a larger constant, if necessary,
we get that the formula above holds for all $t>0$.  This proves Proposition 3.5.\\ \\

\noindent {\bf Proof of Proposition 3.6.} Fix $f\in G\sm\{ 1_M\}$. Let $M_1,... M_k$ the the components of
the fixed-point set of $f$. Then the $M_i$'s are closed connected totally geodesic submanifolds of $M$ and
$dim\, M_i<n$. Define $R^f(t)$ as the set of all
non-fixed oriented simple closed geodesics $\gamma$ in $M$, of length $\leq t$,
with no translations, for which $f|_{\gamma}$ is not trivial.
As mentioned in item 5 above, for such geodesics we have $G_\gamma=\{ 1_M, f\}$ and $f|_\gamma$ reverses orientation.
By definition we have that $P^f(t)$ is the cardinality of $R^f(t)$.\\

For each $\gamma\in R^f(t)$, we have that the fixed-point set of $f|_{\gamma}$ has exactly two points. 
Order this set (arbitrarily) and denote the first point by $p_\gamma$ and
the second one by $q_\gamma$. Define
$$
R_{i,j}(t)\,\, =\,\, \Big\{ \,\,\gamma\in R^f(t)\,\, :\,\, p_\gamma\in M_i,\,\, q_\gamma\in M_j\,\,\Big\}
$$

\noindent and let $P_{i,j}(t)$ be the cardinality of this set. Then $R^f(t)=\bigcup_{i,j}R_{i,j}(t)$
(disjoint union) and $P^f(t)=\sum_{i,j}P_{i,j}(t)$. Therefore it is enough to prove that 
for each $(i,\, j)$ there is a constant $c>0$ such that 
$$
 P_{i,j}(t)\,\,\, \leq\,\,\, c\,\,\frac{e^{\, (2/3)ht}}{t}
$$
\\

Parametrize each $\gamma\in R_{i,j}(t)$ by a constant speed parametrization $\gamma: [0,2]\ra M$ (we use the same letter)
with $\gamma(0)=\gamma(2)=p_\gamma$ and $\gamma(1)=q_\gamma$. Let $p:\Ha^n\ra M$ be the covering projection and identify
the universal cover $\widetilde{M}_i$ of $M_i$ with one of the components of $p^{-1}(M)$. \\

Let $s>0$ be such that $2s$ is less than the widths of the normal tubular neighborhoods of $M_i$ and $M_j$
and less that the injectivity radius of $M$.
Let $C\sbs \widetilde{M}_i$
be a compact set such that $p(C)=M_i$. Let $D$ be the diameter of $C$.
For each $\gamma\in R_{i,j}(t)$ chose a lifting $\widetilde{\gamma}: [0,1]\ra \Ha^n$ of $\gamma|_{[0,1]}: [0,1]\ra M$,
with $\widetilde{\gamma}(0)\in C$. Denote by $B_\gamma$ the closed ball centered at $\widetilde{\gamma}(1)$,
of radius $s$.\\

\noindent {\bf Claim.} {\it If $\gamma\neq \gamma'$ then $B_{\gamma}$ and $B_{\gamma'}$ are disjoint.}\\

\noindent {\bf Proof of the Claim.} Since $\gamma(1)=q_\gamma$, $\gamma'(1)=q_{\gamma'}\in M_j$, 
we have that $\widetilde{\gamma}(1)\in N$ and $\widetilde{\gamma'}(1)\in N'$, where $N$ and $N'$ are
components of $p^{-1}(M_j)$. Note that both $\widetilde{\gamma}$ and $\widetilde{\gamma}'$ are orthogonal
to $\widetilde{M_i}$ (at $\widetilde{\gamma}(0)$ and $\widetilde{\gamma}'(0)$ respectively). Hence, 
if $N=N'$, then $\widetilde{\gamma}$ and $\widetilde{\gamma}$
would both also be orthogonal to $N$, and it follows that $\gamma=\gamma'$. Hence $N\neq N'$. The Claim
follows now from the choice of the number $s>0$. This proves the Claim.\\

Choose a point $x_0\in C$. Note that, for each $\gamma\in R_{i,j}(t)$, the length of $\widetilde{\gamma}|_{[0,1]}$ is 
$\leq t/2$.
Then, for such $\gamma$'s we have
$$
B_{\gamma}\sbs {\mbox{\bf B}}_{\Ha^n}\Big( x_0, (t/2)+s+D\Big)
$$

\noindent Hence 

$$
\bigcup_{\gamma\in R_{i,j}(t)} B_{\gamma}\sbs {\mbox{\bf B}}_{\Ha^n}\Big( x_0, (t/2)+s+D\Big)
$$

\noindent This, together with the Claim imply

$$
\nu\, P_{i,j}(t)=\sum_{\gamma\in R_{i,j}(t)} {\mbox{{\bf Vol}}}\Big(B_{\gamma}\Big)\,\,\, \leq\,\,\, 
{\mbox{{\bf Vol}}}\Bigg({\mbox{\bf B}}_{\Ha^n}\Big( x_0, (t/2)+s+D\Big)\Bigg)
$$

\noindent where $\nu>0$ is the volume of a ball of radius $s$ in $\Ha^n$. Hence, for $t$ large we have

$$
 P_{i,j}(t)\,\,\, \leq\,\,\,\frac{1}{\nu}\,\,
{\mbox{{\bf Vol}}}\Bigg({\mbox{\bf B}}_{\Ha^n}\Big( x_0, (8/15)t\Big)\Bigg)
$$

Now the rest of the proof is similar to the last part of the proof of Proposition 3.5.
This proves Proposition 3.6.
\vspace{1.5in}

\noindent {\bf \Large  Section 4. Proof of Lemma 2.1.}\\

The space $\met(M)$ of all Riemannian metrics -with the smooth topology- is an open set in the separable Fr\'echet
space of all symmetric 2-tensors on $M$. Hence $\met(M)$ is a separable ANR
and, since $\mo(M)$ is open in $\met(M)$, we have that $\mo(M)$ is also a separable
ANR, see \cite{Palais} (here by ANR we mean ``ANR for
metrizable spaces"). By Ebin's Slice Theorem \cite{Ebin} together with results of Borel \cite{Borel} and Conner and Raymond
\cite{CR} we have that $\met(M)\ra {\cal{T}}(M)$ (and thus $\mo(M)\ra\To(M)$) is a locally trivial bundle (see also Lemma 1.1 of \cite{FO5}). Therefore $\To(M)$ is locally an ANR, hence it is itself an ANR, provided it is metrizable. 
Since, by a result of W.H.C. Whitehead \cite{Palais},  every ANR has the homotopy type of a CW complex, Lemma 2.1 follows,
provided we prove that ${\cal{T}}(M)$ (and hence $\To(M)$) is metrizable.\\

\noindent {\bf Remark.} In fact, since they are separable, both $\To(M)$ and $\cT(M)$ have the homotopy type
of a countable CW-complex \cite{Milnor}.\\

\noindent {\bf Proposition 4.1} {\it ${\cal{T}}(M)$  is metrizable.}\\

\noindent {\bf Proof.} Write  ${\bf M}=\met(M)$ and ${\bf D}=DIFF_0(M)$ (both with the smooth topology). Also  ${\bf M}^{C^k}$
is the space of $C^k$-metrics on $M$, with the $C^k$-topology and   ${\bf M}^s$ the Sobolev $s$-completion of ${\bf M}$
(see \cite{Ebin}). Also ${\bf D}^s$ and ${\bf D}^{C^k}$ are defined similarly.
We have the following facts (see \cite{Ebin}):
\begin{enumerate}
\item[(i)] the inclusions ${\bf M}\ra {\bf M}^s$, ${\bf D}\ra {\bf D}^s$ are continuous, and the restriction of the action of
${\bf D}^s$ on ${\bf M}^s$ coincides with the action of ${\bf D}$ on ${\bf M}$.

\item[(ii)] For $s>k+\frac{n}{2}$ we have ${\bf M}^s\sbs {\bf M}^{C^k}$, and the inclusion is continuous.

\item[(iii)] There is a ${\bf D}^s$-invariant metric $d^s$ on ${\bf M}^s$. (This metric generates the topology on ${\bf M}^s$).


\end{enumerate}

\noindent We assume that $d^s\leq 1$ (just replace $d^s$ by $\frac{d^s}{1+d^s}$.) We define the metric $d$ on 
${\bf M}$:

$$
d=\sum_{s\geq 0}\frac{1}{2^s}d^s
$$

\noindent By (i), each $d^s$ is continuous on ${\bf M}\times {\bf M}$, hence $d$ is also continuous on ${\bf M}\times {\bf M}$.
And by (i) and (iii), $d$ is ${\bf D}$-invariant.
We prove that $d$ generates the topology of ${\bf M}$.\\

First, since $d$ is continuous we have that the smooth-topology is finer that the $d$-topology. We have to prove that
the $d$-topology contains the $C^k$-topology (on ${\bf M}$), for every $k$. Let $U$ be an  $C^k$-open set in ${\bf M}$ and $g\in U$.
Then $U=U'\cap {\bf M}$, for some open set $U'$ in ${\bf M}^{C^k}$. 
We prove that there is a $d$-ball centered at $g$ and contained in $U$.
By (ii), for $s$ large we have that
$U'\cap {\bf M}^s$ is open in ${\bf M}^s$. Therefore there is $r>0$ such that:
$$
B^{d_s}(g,r)\,\,\, \sbs\,\,\, U'\cap {\bf M}^s
$$

\noindent where $B^{d_s}(g,r)$ is the $d^s$-ball centered at $g$ of radius $r$. Intersecting with ${\bf M}$ we get:
$$
B^{d_s}(g,r)\cap  {\bf M}\,\,\,\sbs\,\,\, U'\cap {\bf M}^s\cap {\bf M}\,\,\,=\,\,\,U
$$

\noindent But, since $d^s\leq 2^s d$ we have that:
$$
B^{d}(g,\frac{r}{2^s})\,\,\,=\,\,\,B^{d}(g,\frac{r}{2^s})\cap  {\bf M}\,\,\,\sbs\,\,\, B^{d_s}(g,r)\cap  {\bf M}\,\,\, \sbs\,\,\, U
$$

\noindent Therefore $d$ generates the topology on ${\bf M}$.
And $d$ induces a metric $d'$ on $\cT(M)={\bf M}/{\bf D}$, by defining  $d'(x,y)= d(\,\xi^{-1}(x), \xi^{-1}(y)\,)=inf \{ d(g,h)\,
\xi(g)=x,\, \xi(h)=y\}$. It is not hard to show that $d'$ is a metric that induces the quotient topology on $\cT(M)={\bf M}/{\bf D}$. 
This proves the Proposition.

\vspace{1.5in}

\noindent {\bf \Large  Section 5.  Proof of Proposition 2.3.}\\

First recall that $P(N)$ is the space of
topological pseudoisotopies of the closed manifold $N$, that is, the space of all homeomorphisms
$N\times I\ra N\times I$, that are the identity on $N\times \{ 0\}$.
The ``take top" map $\tau:P(X)\ra TOP(X)$ is $\tau(\phi)(x)=\phi(x,1)$.\\

In this Section we will need some objects defined in \cite{FO7}. We will denote the given hyperbolic metric
on $M$ by $g_0$. Let $\epsilon >0$. In Section 2.1 of \cite{FO7} the space $\met_\gamma^\epsilon(M)$ is defined
as the subspace of $\mo(M)$ of metrics $g$ such that: (1) the $g$-geodesic homotopic to $\gamma$ is 
$\gamma$ itself, and (2) the $g$-normal-bundle (to $\gamma$) and $g$-exponential map (at points in
$\gamma$) in an $\epsilon$-neighborhood of $\gamma$ coincide with the ones given by $g_0$.
(For a precise definition see Section 2.1 of \cite{FO7}).\\

Also in \cite{FO7}, a map $\Lambda_\gamma^\epsilon :\met_\gamma^\epsilon(M)\ra P(\bS^1\times\bS^{n-2})$
is defined in, roughly, the following way. Let $Q$ be the covering space of $M$ corresponding to the infinite
cyclic group generated by $\gamma$. We will denote the lifting of any metric on $M$ to $Q$ with the same letter.
Identify $Q$ with $\bS^1\times\R^{n-1}$ such that: (1) $t\mapsto(cos(\frac{2\pi t}{\ell}),sin(\frac{2\pi t}{\ell}),0)
\in\bS^1\times\{ 0\}$
is the (unique) lifting of $\gamma$, (2) $t\mapsto (z,tv)$, $(z,v)\in\bS^1\times\bS^{n-2}$, are speed-one $g_0$-geodesic
rays $g_0$-orthogonal to $\bS^1\times\{ 0\}$. Then the complement of the $\epsilon$-neighborhood of $\bS^1\times\{ 0\}$
(with respect to $g_0$) can be identified with  $\bS^1\times \bS^{n-2}\times [\epsilon, \infty)$. For 
$g\in  \met_\gamma^\epsilon(M)$, $\Lambda_\gamma^\epsilon(g)$ is defined as $exp^g$, i.e, the orthogonal 
(to $\bS^1\times\{ 0\}$) exponential with respect to $g$. Note that $\Lambda_\gamma^\epsilon(g)$
restricted to $\bS^1\times \bS^{n-2}\times \{ \epsilon\}$ is the identity. We can extend, using asymptotics,
$\Lambda_\gamma^\epsilon(g)$ to the ``space at infinity" $\bS^1\times \bS^{n-2}\times \{ \infty\}$
obtaining $\Lambda_\gamma^\epsilon(g):\bS^1\times \bS^{n-2}\times [ \epsilon, \infty]\ra\bS^1\times 
\bS^{n-2}\times [ \epsilon, \infty]$. And, identifying $ [ \epsilon, \infty]$ with $[0,1]$ we can write
$\Lambda_\gamma^\epsilon(g)\in P(\bS^1\times \bS^{n-2})$.
(For more details see Section 2.1 of \cite{FO7}.)\\

Thus, we have the following sequence:
$$\met_\gamma^\epsilon(M)\stackrel{\Lambda_\gamma^\epsilon(g)}{\longrightarrow}
P\Big(\bS^1\times \bS^{n-2}\Big)\stackrel{\tau}{\longrightarrow}    TOP\Big(\bS^1\times \bS^{n-2}\Big)$$
\noindent And we want to relate this sequence to the meaningful sequence defined earlier:
$$\mo(M)\stackrel{\xi}{\longrightarrow}\To(M)\stackrel{\delta}{\longrightarrow}TOP\Big( \bS^1\times \bS^{n-2} \Big)//L $$
\noindent where, recall, $\delta$ is the composition $ \To(M)\stackrel{\eta}{\ra}\To(M)^\bullet
\stackrel{\Delta_\gamma}{\ra}TOP( \bS^1\times \bS^{n-2})//L$. The relationship between these two sequences is
given by the following homotopy commutative diagram, where the two sequences embed:\\

$$\begin{array}{ccccc}
&&\met_\gamma^\epsilon(M)&\stackrel{\tau\Lambda_\gamma^\epsilon(g)}{\longrightarrow}&
TOP\Big(\bS^1\times \bS^{n-2}\Big)\\ \\
&\swarrow&\uparrow&&\uparrow\\ \\
\mo(M)&& \met_\gamma^\epsilon(M)^\bullet&\stackrel{\Big(\tau\Lambda_\gamma^\epsilon(g)\Big)^\bullet}{\longrightarrow}&
TOP\Big(\bS^1\times \bS^{n-2}\Big)^\bullet\\ \\
\xi\downarrow&&\xi^\bullet\downarrow\,\,\,\,\,\,\,\,\,&&\downarrow\\ \\
\To(M)&\stackrel{\eta}{\longrightarrow}&\To(M)^\bullet&\stackrel{\Delta_\gamma}{\longrightarrow}&TOP( \bS^1\times \bS^{n-2})//L
\end{array}$$\\

\hspace{1.8in}{\it Diagram 5.1.}\\ \\\\

\noindent This diagram is a version of the diagram at the end of Section 2 of \cite{FO7}.
Diagram 5.1 will be used later. We want now to see what is the effect of applying $\Lambda_\gamma^\epsilon$
to $f\alpha$, where $f\in G$ and $\alpha$ is as in Section 2 above.\\

Recall that we are assuming the normal bundle of $\gamma$ to be trivial, and we also assume that a given
trivialization of this bundle is given (as in  the beginning of Section 1.2 of \cite{FO7}).
So, we can identify, via this trivialization, the unit sphere bundle of $\gamma$ with $\bS^1\times\bS^{n-2}$. 
Let $f\in G_\gamma$. We denote by $Df:\bS^1\times\bS^{n-2}\ra\bS^1\times\bS^{n-2}$ the derivative of $f$
restricted to the unit sphere bundle. Write $Ef=Df\times 1_I$.\\

\noindent {\bf Remark.} Note that $Df$ has the form $Df(z, v)= (f(z), df_z(v))$, $(z,v)\in \bS^1\times\bS^{n-2}$.
(Here $df_z$ is the derivative of $f$, at $z$.) 
And, since $f$ is an isometry, we have
$Df(z,v)=(Cz, B(z)v)$, for some $C\in O(2)$, and $B:\bS^1\ra O(n-1)$. And if $f\in G_\gamma^+$, then
$Df\in L$. \\

\noindent {\bf Lemma 5.2.}  {\it For $f\in G_\gamma$ and $g\in \met_\gamma^\epsilon(M)$ we have
$\Lambda_\gamma^\epsilon (fg)\,\simeq\, Ef[\Lambda_\gamma^\epsilon (g)] Ef^{-1}$.}\\

\noindent {\bf Proof.} Denote the lifting of $f$ to $Q$ also by $f$. Thus $f:(Q,g_0)\ra (Q,g_0)$ is an isometry,  and $f$ sends
$g_0$-geodesics to $g_0$-geodesics. Because of the identification made above between $Q$ and
$\bS^1\times\R^{n-1}$, we have that $f(z,v)=(f(z), f_*(v))$, for $(z,v)\in \bS^1\times\R^{n-1}$.
Therefore $Df$ is just the restriction of $f$ to $\bS^1\times\bS^{n-2}$. Thus, considering the identifications already made,
we can write $f|_{\bS^1\times \bS^{n-2}\times [\epsilon, \infty]}=Ef$.\\

Now, let $g\in\met_\gamma^\epsilon(M)$ and, as before, denote also by $g$ its lifting to $Q= \bS^1\times\R^{n-1}$.
Since $f:(Q,g)\ra (Q,fg)$ is an isometry,  $f$ sends
$g$-geodesics to $fg$-geodesics. Hence, for $(z,v)\in \bS^1\times\bS^{n-2}$, $t\geq 0$, we have
$$f\Bigg( exp\,^g\Big( f^{-1}(z), tf_*^{-1}(v)\Big)\Bigg)=exp\,^{fg}(z,tv)$$

\noindent Therefore $f\circ exp \,^g\circ f^{-1}=exp\,^{fg}$, and the Lemma follows.\\

\noindent {\bf Lemma 5.3.}  {\it If $f\in G_\gamma^+$, then $Df$ is isotopic to the identity, provided $o(\gamma)$ is odd.}\\

\noindent {\bf Proof.} By the Remark above, $Df$ has the form
$Df(z,v)=(Cz, B(z)v)$, for some $C\in SO(2)$, and $B:\bS^1\ra SO(n-1)$. Hence $Df\in SO(2)\times\Omega SO(n-1)\sbs TOP(\bS^1\times\bS^{n-2})$. Note that $\pi_0(SO(2)\times\Omega SO(n-1))=\pi_1(SO(n-1))=\Z_2$, $n>3$ (and $\Z$,
for $n=3$). Therefore, since $G_\gamma^+\ra
SO(2)\times\Omega SO(n-1)\sbs TOP(\bS^1\times\bS^{n-2})$, $f\mapsto Df$, is a group homomorphism and 
$|G_\gamma^+|$ is odd, $Df$ belongs to the identity component of $SO(2)\times\Omega SO(n-1)$. Hence
$Df$ is isotopic to the identity. This proves the Lemma.\\

\noindent {\bf Lemma 5.4.}  {\it We can choose $\alpha$ so that it satisfies the following
two conditions:}
\begin{enumerate}
\item[{\bf (a.)}] {\it If $f\notin G_\gamma$ then  $\Lambda_\gamma^\epsilon f\alpha$ is null-homotopic.}

\item[{\bf (b.)}] {\it If $f\in G_\gamma^+$ then $\Lambda_\gamma^\epsilon f\alpha\,\simeq\, \Lambda_\gamma^\epsilon \alpha$.}
\end{enumerate}

\noindent {\bf Proof.} 
Part (b) follows from Lemmas 5.2 and 5.3.
We now prove part (a). For this we use the notation of Section 6 of \cite{FO7}.
Let $\alpha:S^k\ra\met_\gamma^\epsilon(M)$ be as constructed in there, and let $f\notin G_\gamma$.
Then the lifting $g_u$ of the metric  $f\alpha (u)$, $u\in\bS^{k}$, to $Q$ almost coincides with the metric
$\bar{\omega}_u$ that appears in the proof of Claim 2, in Section 6 of \cite{FO7}. These metrics
differ only on the set $R$, where $g_u$ is equal to the given hyperbolic metric. Then deformations
similar to the two deformations in the proof of Claim 2, in Section 6 of \cite{FO7} and the argument used
in Claim 3 prove part (a) of the Lemma.\\

 Proposition 2.3. now follows from Lemma 5.4, Diagram 5.1 and the fact that the upward arrows in
Diagram 5.1  are weak homotopy equivalences.

F.T. Farrell

SUNY, Binghamton, N.Y., 13902, U.S.A.\\

P. Ontaneda

SUNY, Binghamton, N.Y., 13902, U.S.A.

\end{document}